\documentclass{article}

\usepackage{arxiv}

\usepackage[utf8]{inputenc} 
\usepackage[T1]{fontenc}    
\usepackage{hyperref}       
\usepackage{url}            
\usepackage{booktabs}       
\usepackage{amsfonts}       
\usepackage{nicefrac}       
\usepackage{microtype}      
\usepackage{lipsum}
\usepackage{amsmath}
\usepackage{amssymb}
\usepackage{cleveref}
\usepackage{graphicx}

\newtheorem{remark}{Remark}

\title{Learning on dynamic statistical manifolds}

\author{
  Francesca Boso \\
  Energy Resources Engineering\\
  Stanford University\\
  Stanford, CA 94305 \\
  \texttt{fboso@stanford.edu} \\
   \And
 Daniel M. Tartakovsky \\
  Energy Resources Engineering\\
  Stanford University\\
  Stanford, CA 94305 \\
  \texttt{tartakovsky@stanford.edu} \\
}

\begin{document}
\maketitle

\begin{abstract}
Hyperbolic balance laws with uncertain (random) parameters and inputs are ubiquitous in science and engineering. Quantification of uncertainty in predictions derived from such laws, and reduction of predictive uncertainty via data assimilation, remain an open challenge. That is due to nonlinearity of governing equations, whose solutions are highly non-Gaussian and often discontinuous. To ameliorate these issues in a computationally efficient way, we use the method of distributions, which here takes the form of a deterministic equation for spatiotemporal evolution of the cumulative distribution function (CDF) of the random system state, as a means of forward uncertainty propagation. Uncertainty reduction is achieved by recasting the standard loss function, i.e., discrepancy between observations and model predictions, in distributional terms. This step exploits the equivalence between minimization of the square error discrepancy and the Kullback-Leibler divergence. The loss function is regularized by adding a Lagrangian constraint enforcing fulfillment of the CDF equation. Minimization is performed sequentially, progressively updating the parameters of the CDF equation as more measurements are assimilated.
\end{abstract}


\section{Introduction}

Robust and efficient quantification of parametric uncert- ainty in hyperbolic balance and conservations laws is hampered by their nonlinearity and solution structure, which typically posses sharp gradients and often exhibits shocks and/or discontinuities. Many uncertainty quantification techniques (e.g., stochastic finite elements and stochastic collocation), which can be orders of magnitude faster than standard Monte Carlo simulations (MCS) when applied to elliptic and parabolic equations, often underperform on hyperbolic problems.

The method of distributions (MD)~\cite{tartakovsky2017method} is an uncertainty quantification technique that is tailor-made for hyperbolic problems with random coefficients and inputs. Its goal is to derive a deterministic partial-differential equation (PDE) for either the probability density function (PDF) or the cumulative distribution function (CDF) of the model output. In the presence of multiplicative noise introduced, e.g., by random parameter fields, MD requires a closure approximation, which is derived either via perturbation expansions or by resorting to phenomenology~\cite{boso-2014-cumulative,boso-2016-method,alawadhi-2018-method}. The method does not rely on a finite-term representation (e.g., via a truncated Karhunen-Lo\`eve expansion) of random parameter fields and, hence, does not suffer from the so-called ``curse of dimensionality'' {\cite{ghanem_2015,tartakovsky2017method}}; its computational cost is independent of the correlation length of an input parameter~\cite{venturi2013exact} and can be orders of magnitude lower than that of MCS~\cite{Yang2019Method,boso-2014-cumulative,alawadhi-2018-method}, and its accuracy increases as the correlation length decreases~\cite{tartakovsky2009probability, tartakovsky2017method}.

While MD enables one to quantify predictive uncertainty in hyperbolic models, assimilation of observations into probabilistic model predictions, e.g., by means of Bayes' rule, facilitates reduction of this uncertainty. 
%
%
Within this framework, the model provides a link between observed quantities and the estimates of the state, filtered through an observational map~\cite{wikle2007bayesian}. 
Direct application of Bayes' rule is often impractical because of the high dimensionality of a joint PDF of system states, and because of complex relations between parameters and states and between states and observations~\cite[sec. 10.2]{evensen2009data}. For these reasons, a plethora of approximation techniques have been proposed. Some of these, e.g., maximum likelihood estimation (MLE)~\cite{myung-2003-tutorial} and maximum a posteriori estimation (MAP)~\cite{cousineau-2013-improving}, aim to identify the mode of a posterior distribution, which can be inadequate if the latter is highly non-Gaussian (e.g., multimodal), as is typical of nonlinear models. Ensemble Kalman filters~\cite{katzfuss-2016-understanding} allow one to handle nonlinear PDEs but assume that their solutions are Gaussian. Other methods, e.g., Markov chain Monte Carlo (MCMC)~\cite{brooks-2011-handbook} and particle filters~\cite{speekenbrink-2016-tutorial}, aim at sampling from the posterior directly and obviate the need for the Gaussianity and linearity assumptions. Like direct Bayesian updating, the methods of this class are computationally expensive because they rely on multiple forward solves of PDEs with uncertain (random) coefficients and/or auxiliary functions. Our goal is to eliminate this step by replacing it with MD.

Variational formulation recasts some of the methods described above (MLE, MAP, analysis step in EnKF) as a minimization problem in which a cost (loss) function contains the average distance between measurements and a model's predictions; parameter estimation is then accomplished by minimizing this loss function with respect to the model's parameters (and their statistical moments). This variational formulation belongs to a broader class of optimization methods, sometimes termed Variational Inference (VI)~\cite{blei2017variational}, that approximate Bayesian posterior densities by imposing closeness (in the Kullback-Leibler divergence sense) to the target density. Key innovations of our method are to reformulate the loss function in distributional terms using a different discrepancy metric and to confine both the prior and the posterior distributions to a dynamic statistical manifold defined by a deterministic CDF equation. 
Minimization is done with respect to variables used to parameterize the closure terms in the CDF equation; these variables are, in turn, expressed in terms of the statistical properties of the uncertain parameters and/or auxiliary functions of the original model. 

Resulting PDE-constrained optimization problems can be solved with several techniques~\cite{herzog-2010-algorithms}. We employ a machine learning approach~\cite{raissi2017machine, zhu2019physics, zhang2019quantifying}, which approximates a PDE's solution with a neural network whose coefficients are obtained by minimizing the resulting residual. 
This component of our algorithm places it in the burgeoning field variously known as physics-informed machine learning or data-aware modeling. Its goal is to overcome the scarcity of experimental data inherent in many physical systems by fusing physical constraints and observations. It is worthwhile emphasizing though that optimization techniques other than the one mentioned above can be used in our Bayesian data assimilation algorithm.

In \cref{sec:forecast}, we formulate a data assimilation problem for hyperbolic PDEs with uncertain parameters and/or auxiliary functions, and introduce  MD as a \textit{forecast} step in Bayesian updating. Section~\ref{sec:analysis} contains a novel \textit{analysis} step, in which MD is used as a constraint to reduce parametric uncertainty{; technical details are provided in \cref{app:MD}}. We refer to this combination of forecast and analysis as the data-aware method of distributions (DA-MD). In \cref{sec:problemformulation}, we test our approach on a linear inhomogeneous hyperbolic equation; this setting admits both exact and approximate Bayesian updates of the random parameters (either spatially uniform or variable) and, hence, enables us to verify the method's accuracy. {Finally, in \cref{sec:conclusion}, we summarize the main findings and discuss future directions.}

\section{Forecast: Method of Distributions}
\label{sec:forecast}

While the data assimilation approach introduced here is applicable to other problems, we formulate it in \cref{sec:forecast}\ref{sec:formulation} for hyperbolic PDEs with uncertain (random) parameters and/or auxiliary functions. This setting simplifies the derivation of a deterministic CDF equation used in \cref{sec:forecast}\ref{sec:MD} as the forecast step in Bayesian data assimilation.

\subsection{Problem Formulation}
\label{sec:formulation}

We consider a {smooth} state variable $u(\mathbf x,t): \Omega \times \mathbb R^+ \rightarrow \mathbb R$, whose dynamics is governed by a nonlinear hyperbolic PDE
\begin{subequations}\label{eq:conslaw}
\begin{align}
  \frac{\partial u}{\partial t} + \nabla \cdot \mathbf q(u; \boldsymbol \theta_q) = r(u; \boldsymbol \theta_r), \qquad \mathbf x \in \Omega, \quad t>0. 
\end{align}
This equation is subject to the initial condition
\begin{align}
 u(\mathbf x,t=0) = u_0(\mathbf x), \qquad \mathbf x \in \Omega 
\end{align}
\end{subequations}
and, if the $d$-dimensional domain $\Omega \subset \mathbb R^d$ is bounded, to appropriate boundary conditions along the domain boundary $\partial \Omega$. The flux, $\mathbf q (u) : \mathbb R \rightarrow \mathbb R^d$, and the source term, $r(u) : \mathbb R \rightarrow \mathbb R$, are parameterized by $\boldsymbol \theta_q$ and $\boldsymbol \theta_r$, respectively. These real-valued parameters can either be constant or vary in space ($\mathbf x$) and time ($t$). The functions $\mathbf q(u)$ and $r(u)$ are either linear or nonlinear, as long as the solution of~\eqref{eq:conslaw} does not develop shocks.\footnote{The presence of shocks and discontinuities complicates the derivation of CDF equations~\cite{alawadhi-2018-method, wang2013cdf, boso-2020-data}, obfuscating our focus on data assimilation.} For example, $u(\mathbf x,t)$ is the concentration of a reactive solute advected by a flow velocity $\mathbf v (\mathbf x)$, while undergoing chemical transformations; in this setting, $\mathbf q(u) = \mathbf v (\mathbf x) u$ is the advective flux parameterized by $\mathbf v(\mathbf x)$, and $r(u)$ represents a chemical reaction parameterized by a reaction rate constant $k$.

Incomplete or noisy measurements of the parameters $\boldsymbol \theta = \{ \boldsymbol \theta_q, \boldsymbol \theta_r \}$ render them uncertain; this uncertainty is quantified by treating $\boldsymbol \theta$ as random fields and random variables. Additionally, auxiliary functions, such as the initial state $u_0(\mathbf x)$ and boundary functions, are uncertain/random. In the following, $\tilde{\boldsymbol \theta}$ denotes the complete set of random inputs, comprised of both $\boldsymbol \theta$ and auxiliary functions. This randomness renders, $u(\mathbf x,t)$,  a solution of~\eqref{eq:conslaw}, random as well. {Rather than computing low statistical moments of $u(\mathbf x,t)$ (e.g., its ensemble mean $\bar u(\mathbf x,t)$ and standard deviation $\sigma_u(\mathbf x,t)$ that are commonly used to obtain an unbiased estimator of a system's dynamics and to quantify the corresponding predictive uncertainty, respectively), our goal is to compute} its one-point CDF $F_u(U;\mathbf x,t) \equiv \mathbb P[u(\mathbf x,t) \le U]$ {where $U \in \Omega_U \subseteq \mathbb R$. The value space for the random variable $u(\mathbf x,t)$, $\Omega_U = [U_{\text{min}}, U_{\text{max}}]$, identifies the support of the CDF $F_u(U;\cdot)$. The latter can be either infinite ($\Omega_U = \mathbb R$, with $U_{\text{min}} = - \infty$ and $U_{\text{max}} = + \infty$) or finite ($U_{\text{min}}, U_{\text{max}} \in \mathbb R$ such that $U_{\text{min}} < U_{\text{max}}$).} 

The model~\eqref{eq:conslaw} is supplemented with $N_\text{meas}$ measurements of the state variable $u(\mathbf x,t)$ collected at selected space-time points $(\mathbf x, t)_m$ with $m = 1,\cdots, N_\text{meas}$. These data, $\mathbf d_{1: N_\text{meas} } = \{d_1, \cdots, d_{N_\text{meas}} \}$, are assumed to differ from the corresponding model predictions $u[(\mathbf x,t)_m]$ by a random measurement error $\varepsilon_m$,
\begin{align}\label{eq:data}
d_m = u[(\mathbf x,t)_m] + \varepsilon_m, \qquad m = 1,\cdots, N_\text{meas}.
\end{align}
The measurement errors are assumed to have zero mean, $\mathbb E[\varepsilon_m] = 0$, and to be mutually {uncorrelated}, $\mathbb E[\varepsilon_m \varepsilon_n] = 0$ for all $m \neq n$. A complete probabilistic description of the data is encapsulated in the PDF 
$f_L(d_m | u[(\mathbf x,t)_m] = U)$, which is also known as likelihood function. In the absence of measurement errors, the observational PDF is given by the Dirac distribution $\delta(\cdot)$, i.e., $f_L(d_m | u[(\mathbf x,t)_m] = U) = \delta (U-d_m)$. 

\subsection{CDF Equation}
\label{sec:MD}

Direct numerical computation of the CDF $F_u(U;\mathbf x,t)$, e.g., via Monte Carlo simulations of~\eqref{eq:conslaw}, is computationally expensive. Instead,  we use MD to derive a $(d+1)$-dimensional linear PDE for $F_u$ (see \cref{app:MD} for details),
\begin{align}\label{eq:CDF}
\frac{\partial F_u}{\partial t} + \boldsymbol {\mathcal Q} (U;\mathbf x,t) \cdot \widetilde \nabla F_u  = \widetilde \nabla \cdot \left[ \boldsymbol {\mathcal D} (U;\mathbf x,t) \widetilde \nabla F_u \right], \qquad (\mathbf x,U) \in \widetilde \Omega, \quad t>0.
\end{align}
This deterministic PDE is defined in the augmented space $\tilde \Omega = \Omega \cup \Omega_U$.This equation is subject to initial and boundary conditions that reflect uncertainty in the initial and boundary conditions for the original problem~\eqref{eq:conslaw}. Additional boundary conditions are defined for $\partial \Omega_U$, $F_u{(U_\text{min};\cdot)} = 0$ and $F_u{(U_\text{max};\cdot)} = 1$; they stem from the definition of a CDF. 

In general, derivation of~\eqref{eq:CDF} requires a closure approximation, such as the perturbation expansion used in \cref{app:MD}. 
Notable exceptions of practical significance include a scenario of random inputs (initial and boundary conditions) but deterministic parameters $\boldsymbol \theta$\footnote{When both the inputs and parameters are deterministic, 
the strategy of transforming a nonlinear $d$-dimensional hyperbolic PDE into its linear $(d+1)$-dimensional counterpart is referred to as kinetic formulation of a hyperbolic conservation law~\cite{perthame-2002-kinetic}. }; in this case~\eqref{eq:CDF} is exact and its coefficients are given by (\cref{app:MD})
\begin{align}\label{eq:exactclosures_deterministicdynamics}
    \boldsymbol {\mathcal Q} (U;\mathbf x,t) = \{ \dot{\mathbf q}(U; \boldsymbol \theta_q), r(U;\mathbf x,t)\}, \qquad \boldsymbol {\mathcal D} (U;\mathbf x,t) = \mathbf 0,
\end{align}
{where $\dot{\mathbf q}(U) = \text d \mathbf q(U) / \text d U$.}
When the model parameters $\boldsymbol \theta$ are random, i.e., when the CDF equation~\eqref{eq:CDF} in inexact, the coefficients $\boldsymbol{\mathcal Q}$ and $\boldsymbol{\mathcal D}$ depend on a set $\boldsymbol \varphi$ of statistical parameters that characterize the randomness of $\boldsymbol \theta$. This set consists of the shape parameters of PDFs of $\tilde{\boldsymbol \theta}$, i.e., their means, variances, and correlation lengths. Together with $(\mathbf x,t)$ and the statistical characteristics of the random auxiliary functions, these parameters represent the coordinates $\tilde{\boldsymbol \varphi}$ of a manifold of distributions $\mathcal F(F_u)$, whose dynamics is governed by the CDF equation~\eqref{eq:CDF}. {Each point in this finite-dimensional coordinate space $\tilde{\boldsymbol\varphi}$ uniquely identifies a distribution \cite{cafaro2011quantifying}.}

The use of perturbative closures to derive a CDF equation raises several questions about its accuracy and robustness, which have been the subject of previous investigations. First, even though the coefficient of variation (CV) of the model parameters serves as a perturbation parameter, the resulting CDF equations for many applications remain accurate for relatively large values of CV~\cite{boso-2014-cumulative,tartakovsky2009probability, Wang2012Uncertainty}. Second, the coefficients of perturbation-based CDF equations, such as $\boldsymbol{\mathcal Q}$ and $\boldsymbol{\mathcal D}$ in~\eqref{eq:CDF}, depend only on the low-order statistical moments (such as $\boldsymbol \varphi$) of the model parameters, rather than their full PDFs. By using an advection-reaction equation as a test-case, we show in \cref{app:MD} that the resulting CDF equation is distributionally robust, giving consistent predictions of the system state's CDF regardless of whether the model coefficient (spatially varying reaction rate) has a Gaussian, log-normal, or uniform PDF. Third, the accuracy of perturbation-based CDF equations depends on correlation lengths of the model parameters: these equations are often exact for white noise (zero correlation) and become progressively less so as the correlation lengths increase. If the correlation lengths are large, perturbation-based closures can be replaced with truncated Karhunen-Lo\'eve expansions of the random parameter fields, leading to accurate/exact CDF equations~\cite{venturi2013exact}. 

In summary, we use the CDF equation~\eqref{eq:CDF} as an efficient forecasting tool, which propagates parametric uncertainty in space and in time through a physical model. It represents a counterpart of a set of ensemble members or particles in the context of ensemble Kalman filter or particle filter, respectively. {Its accuracy and computational efficiency vis-\`a-vis Monte Carlo simulations have been throughly investigated~\cite{boso-2014-cumulative,alawadhi-2018-method,Yang2019Method}.}


\section{Analysis: Sequential Bayesian Update on Dynamic Manifolds}
\label{sec:analysis}

We use MD as a constraint for the analysis step, during which observations of the system state are used to refine the knowledge of the meta-parameters $\boldsymbol \varphi$. Specifically, our novel analysis step involves minimization of the discrepancy between the ``observational'' CDF $\hat F_u(U; (\mathbf x,t)_m)$ in each measurement location ($m=1,\dots,N_\text{meas}$) and the corresponding ``estimate'' CDF $F_u(U; \boldsymbol \varphi; (\mathbf x,t)_m)$:
\begin{align}\label{eq:loss}
   \boldsymbol \varphi^{(m)} = \underset{{\boldsymbol \varphi}}{\text{argmin} } \, \| \hat F_u(U; (\mathbf x,t)_m) - F_u(U; \boldsymbol \varphi; (\mathbf x,t)_m) \|_{2} \quad \text{subject to} \quad F_u \in \mathcal F,
\end{align}
{where 
\[ \| \hat F_u(U; (\mathbf x,t)_m) - F_u(U; \boldsymbol \varphi; (\mathbf x,t)_m) \|_{2} = \left( \int_{\Omega_U} ( \hat F_u(U; (\mathbf x,t)_m) - F_u(U; \boldsymbol \varphi; (\mathbf x,t)_m) )^2 \text d U \right)^{1/2}.
\]}
The analysis step, i.e., minimization of~\eqref{eq:loss}, is performed sequentially for each observation $m$, so that all the distributions above are uni-variate. Formulation \eqref{eq:loss} is at the core of our data assimilation strategy and requires a thorough explanation.

\begin{remark}

\textit{MD constraint}: The estimate distribution $F_u(U; \boldsymbol \varphi; (\mathbf x,t)_m)$ is a solution of the CDF equation~\eqref{eq:CDF} subject to appropriate initial/boundary conditions. This boundary value problem is parameterized by the set of parameters ${{\boldsymbol \varphi}}$,  over which the discrepancy minimization is performed.  In other words,~\eqref{eq:loss} identifies the parameters of the CDF equation that yield a CDF $F_u$ in the measurement location as close as possible to the observational CDF $\hat F_u$. This implies that the minimization is performed on the manifold of distributions obeying the CDF equation. This observation is further elaborated upon in \cref{sec:analysis}\ref{sec:ig}. Reliance on MD obviates the need for both Gaussianity assumption for the system states and the linearity requirement for the physical model, as long as it is possible to develop a reliable and accurate CDF equation.

\end{remark}

\begin{remark}

\textit{Observational CDFs}:  We construct the observational CDF, 
\[
\hat F_u(U; (\mathbf x,t)_m) = \int_{U_\text{min}}^U \hat f_u(U; (\mathbf x,t)_m) \text d U,
\]
 via Bayesian update of the corresponding PDF $\hat f_u$ at each space-time measurement point $m$:
\begin{align}\label{eq:bayesloc}
    \hat f_u (U; (\mathbf x,t)_m | d_m) \propto f_L( d_m | u[(\mathbf x,t)_m] = U) f_u(U; {\boldsymbol \varphi}^{(m-1)}; (\mathbf x,t)_m).
\end{align}
The PDF $f_u(U; {\boldsymbol \varphi}^{(m-1)}; (\mathbf x,t)_m)$ is computed from a solution of the CDF equation~\eqref{eq:CDF} whose parameters ${{\boldsymbol \varphi}^{(m-1)}}$ are computed in the previous assimilation step. 
This procedure provides a \emph{local} update of the system state's PDF in the sense that it yields no information on the surrounding locations nor on the future time evolution of the state.  

\end{remark}

\begin{remark}

\textit{Sequential update}: The sequential update of the observational PDF $f_L$ allows us to obtain final estimates for the MD parameters ${\boldsymbol \varphi}$ that are conditional on all assimilated observations {\cite{giffin2007updating}}. It is employed both to reduce the dimensionality of the CDFs/PDFs involved and to facilitate real-time update of the estimates as new measurements become available~\cite[p.~101]{evensen2009data}. At each step, or for each data point, $m = 1, \dots, N_\text{meas}$, we follow the following procedure.
\begin{itemize}
        \item For $m=1$, the MD parameters $\boldsymbol{\varphi}^{(0)}$ are initialized to define the prior and to compute~\eqref{eq:bayesloc}. The normalization constant in~\eqref{eq:bayesloc} is obtained by (numerical) integration, $C_1 = \int  f_L(d_1 | U) f_u(U; {\boldsymbol \varphi}^{(0)}; (\mathbf x,t)_1) \text d U$.
        \item For $m>1$, each update~\eqref{eq:bayesloc} accounts for conditioning on all previous measurements up to the current one, $\mathbf d_{1:m}$, such that
        \begin{align}
            \hat f_u (U; (\mathbf x,t)_m | \mathbf d_{1:m}) \propto f_L(\mathbf d_{1:m} | U) f_u (U; {\boldsymbol \varphi}^{(m-1)}; (\mathbf x,t)_m).
        \end{align}
        This step implies that the prior distribution in the current measurement location $m$ obeys the CDF equation~\eqref{eq:CDF}. If observation errors are mutually uncorrelated, then $f_L(\mathbf d_{1:m} | U) = \prod_{i=1}^m f_L(d_i | U)$ and
        \begin{align}\label{eq:postloc}
    \hat f_u (U; (\mathbf x,t)_m | \mathbf d_{1:m})  & \propto \prod_{i=1}^{m-1} f_L (d_i | U) f_L (d_m | U) f_u (U; {\boldsymbol \varphi}^{(m-1)}; (\mathbf x,t)_m) \notag \\
    & \propto f_L(d_m | U) \hat f_u(U; (\mathbf x,t)_m|\mathbf d_{1:m-1}).
\end{align}
	Here, $\hat f_u(U; (\mathbf x,t)_m|\mathbf d_{1:m-1})$ is approximated by a solution of the CDF equation in $(\mathbf x,t)_m$ with parameters $\boldsymbol \varphi^{(m-1)}$ from the previous iterative step. In other words, a solution of the CDF equation~\eqref{eq:CDF} with parameters $\boldsymbol \varphi^{(m-1)}$ serves as prior.
\end{itemize} 
At the end of this sequential assimilation procedure, the CDF equation~\eqref{eq:CDF} with parameters ${\boldsymbol \varphi}^{(N_\text{meas})}$ allows us to predict the future dynamics of the CDF $F_u(U;\cdot)$, i.e., to make a probabilistic forecast.

\end{remark}

\begin{remark} 
\textit{Choice of the discrepancy metric}: Our reliance on the squared $L^2$ norm {(a.k.a. Cramer's distance \cite{g.2018the})}, 
\[
\| F_1(U) - F_2(U)\|_2^2 = \int_{U_\text{min}}^{U_\text{max}} [F_1(U) - F_2(U)]^2 \text d U,
\]
as a measure of discrepancy between any two CDFs, $F_1(U)$ and $F_2(U)$, facilitates numerical minimization of the loss function in~\eqref{eq:loss} with a technique described in \cref{sec:analysis}\ref{sec:minimization} below. We deploy it in place of a commonly used Kullback-Leibler (KL) divergence, 
\[
  D_\text{KL}(F_1, F_2) = \int_{U_\text{min}}^{U_\text{max}} f_1(U) \ln \frac{f_1(U)}{f_2(U)} \text d U,
\]
for the following reasons. 
According to 
Pinsker's inequality~\cite{pinsker-1964-information,Topsoe-2000-some},
$D_\text{KL}[F_1, F_2] \ge (1/2) \| F_1 - F_2\|_1^2$
where $\| \cdot \|_1$ is the $L^1$ norm. Since $\|F_1 - F_2 \|_1 \ge \| F_1 - F_2\|_2$ {\cite[Prop. 1.5]{stein2011functional}}, this yields $D_\text{KL}(F_1, F_2) \ge (1/2) \| F_1 - F_2\|_2^2$.
Since $D_{\text{KL}}(F_1,F_2)$ and $\| F_1 - F_2\|_2$ share the same minimum (for $F_1 \equiv F_2$ both metrics are equal to zero),   a solution of the minimization problem~\eqref{eq:loss} would also minimize the corresponding loss function based on the KL divergence. Moreover, it is advantageous to employ MD in its CDF form, rather than its PDF form, because of the straightforward assignment of the boundary conditions along $\partial \Omega_U$ and smoother solutions.

\end{remark}

\begin{remark} \textit{Relationship to Variational Inference Techniques}: 
Our method aims at approximating posterior densities in a Bayesian sense via a minimization procedure. As such, it connects with VI techniques, which use optimization to identify one joint density---chosen to belong to a specified family of approximate densities---which is close to the target posterior in KL divergence terms~\cite{blei2017variational}. We choose a physics-based family of plausible distributions, which obey the CDF equation parameterized with a finite set of parameters. Constraining distributions to a dynamic manifold allows us to consider sequentially the update of single-point distributions: updated parameters can be used, in combination with the CDF equation, to obtain forecast predictions in different space-time locations. Moreover, it reduces drastically (to one) the dimensionality of the posterior distribution to be updated at each assimilation step.\footnote{{In this regard, we mention the work by \cite{cafaro2011quantifying}, where the reduction in complexity of statistical models is quantified by exploiting relevant embedding constraints specifying geodesic motion on curved statistical manifolds.}}

\end{remark}

\subsection{Loss Function Minimization}\label{sec:minimization}

The PDE-constrained optimization problem~\eqref{eq:loss} can be solved with several techniques~\cite{herzog-2010-algorithms}. If the CDF equation~\eqref{eq:CDF} admits an analytical solution, e.g., if the system parameters $\boldsymbol\theta$ are deterministic and the initial and/or boundary functions are random, $F_u(U;\boldsymbol\varphi)$ can be expressed as a (semi)explicit function of the statistical parameters, $\boldsymbol \varphi_0$ and $\boldsymbol \varphi_b$, characterizing the initial and boundary CDFs $F_0$ and $F_b$, respectively. 
Section~\ref{sec:problemformulation}\ref{sec:randominputs} deals with such a scenario; it serves to verify the reliability of our approach by comparing its performance with that of the standard Bayesian update.


When the CDF equation~\eqref{eq:CDF} has to be solved numerically, 
we follow~\cite{raissi2019physics, zhu2019physics} to approximate its solution, $F_u(U;\tilde{\boldsymbol \varphi})$, with a neural network $F_\text{NN}(U;\tilde{\boldsymbol \varphi})$ whose coefficients (weights and biases) are computed by minimizing the residual
\begin{equation}\label{eq:res}
    R = \frac{\partial F_\text{NN} }{\partial t} + (\boldsymbol{\mathcal Q} - \tilde \nabla \cdot \boldsymbol{\mathcal D} ) \cdot \tilde \nabla F_\text{NN} -  \boldsymbol{\mathcal D} \tilde \Delta F_\text{NN}
\end{equation}
at a set of $N_\text{res}$ points $\{(\mathbf x,t)_r\}_{r=1}^{N_\text{res}}$; the initial and boundary conditions are enforced at a finite set of $N_\text{aux}$ points $\{(U,\mathbf x,t)_r\}_{r=1}^{N_\text{aux}}$.  
The derivatives in~\eqref{eq:res} are computed via automatic differentiation, as implemented in \texttt{TensorFlow}~\cite{tensorflow2015-whitepaper}.
This procedure replaces the PDE-constrained minimization problem~\eqref{eq:loss} with an optimization problem
\begin{align}\label{eq:loss2}
  \boldsymbol \varphi^{(m)} = \underset{{\boldsymbol \varphi}}{\text{argmin} } \big\{ \| \hat F(U;(\mathbf x,t)_m) - F_{\text{NN}}(U;(\mathbf x,t)_m,\boldsymbol{\varphi})\|_{2} + \text{MSE}_R(\boldsymbol \varphi) + \text{MSE}_B(\boldsymbol \varphi) \big\},
\end{align}
where
\begin{align}
    & \text{MSE}_R (\boldsymbol \varphi) = \frac{1}{N_\text{res}} \sum_{r=1}^{N_\text{res}} \| R((\mathbf x,t)_r; \boldsymbol \varphi) \|_{2},  \notag \\
    & \text{MSE}_B(\boldsymbol \varphi) = \frac{1}{N_\text{aux}} \sum_{i=1}^{N_\text{aux}} \| F_{\text{NN}}((U,\mathbf x,t)_i, \boldsymbol \varphi) - F_{\text{inp}}((U,\mathbf x,t)_i)\|_{2}, \notag
\end{align}
where $F_{\text{inp}}$ represents the prescribed CDFs of either the initial state or the boundary functions along $\partial \tilde \Omega$. The NN function approximation via minimization enjoys  convergence guarantees in the chosen $L^2$ norm, e.g., \cite{barron1993universal, bolcskei2019optimal}. A solution of~\eqref{eq:loss2} provides a CDF surrogate (a ``trained'' NN) and the set of optimal parameters~$\boldsymbol \varphi$. The surrogate can then be used to update predictions and for forecast (not pursued here).

\subsection{Information-Geometric Interpretation}
\label{sec:ig}

A family of distributions satisfying the CDF equation~\eqref{eq:CDF} defines a dynamic statistical manifold $\mathcal F[F_u; \tilde{\boldsymbol \varphi}]$. Each point in this space, with coordinates $\tilde{\boldsymbol \varphi} = (\mathbf x,t,\boldsymbol \varphi)$, uniquely identifies a physics-informed CDF $F_u(U; \mathbf x,t)$ of the model's output $u(\mathbf x,t)$ at each space-time point $(\mathbf x,t)$. The manifold $\mathcal F$ is differentiable in all coordinate directions and equipped with a Riemannian metric. The latter takes the form of the Fisher information metric (FIM), a $(d+1+N_\varphi) \times (d+1+N_\varphi)$ matrix whose components are~\cite[p.~33]{amari-2016-information}
\begin{align}
    g_{jk}(\tilde{\boldsymbol \varphi}) = \int \frac{\partial \ln f_u(U;\tilde{\boldsymbol \varphi})}{\partial \tilde \varphi_j} \frac{\partial \ln f_u(U;\tilde{\boldsymbol \varphi})}{\partial \tilde \varphi_k} f_u(U;\tilde{\boldsymbol \varphi}) \text d U, \qquad j,k = 1,\dots,d+1+N_\varphi,
\end{align}
where $N_{\tilde \varphi} = d+1+N_\varphi$ is the number of manifold coordinates, with $N_{\varphi}$ statistical parameters in the CDF equation~\eqref{eq:CDF}.\footnote{{This definition assumes the existence of the PDF $f_u$; for hyperbolic PDEs~\eqref{eq:conslaw} with smooth solutions, it does exist and satisfies a PDF equation corresponding to the CDF equation~\eqref{eq:CDF}~\cite{tartakovsky2017method, venturi2013exact, tartakovsky2009probability}.}} The local curvature of the manifold, $g_{jk}$, represents a Euclidean metric (a distance on the manifold $\mathcal F$) upon an appropriate change of variable. FIM quantifies the differential amount of information between two infinitesimally close points on a manifold; it is formally computed as the second derivative of the KL divergence of distributions $F_u(U;\tilde{\boldsymbol \varphi})$ and $F_u(U;\tilde{\boldsymbol \varphi}')$  with $\tilde{\boldsymbol \varphi}' \rightarrow \tilde{\boldsymbol \varphi}$ \cite{kullback1997information}. 

The significance of FIM and its geometric implications~{\cite{cafaro2020information}} will be explored elsewhere. Here we focus on the calculation of the information gain achieved during each step of the data assimilation process.
Specifically, we express an $m$th analysis step in geometrical terms as a change of the coordinates on the statistical manifold $\mathcal F$, from $\tilde{\boldsymbol \varphi}^{(m-1)}$ to $\tilde{\boldsymbol {\varphi}}^{(m)}$, and quantify the corresponding information gain by $D_{\text{KL}}[F_u(U; \varphi^{(m)}), F_u(U;\varphi^{(m-1)})]$. {This quantity is computed as a post-processing step for comparative analysis.}






\section{Numerical Experiments}\label{sec:problemformulation}

Let us consider a scalar $u(x,t) : \mathbb R^+ \times \mathbb R^+ \rightarrow \mathbb R^+$, whose dynamics satisfies a one-dimensional {dimensionless} advection-reaction equation
\begin{subequations}\label{eq:case0}
\begin{align}
  \frac{\partial u}{\partial t} + \frac{\partial q(u)}{\partial x} = r(x,u), \qquad q \equiv v u, \quad r \equiv - k(x) u; \quad x>0, \quad t>0,
\end{align}
subject to initial and boundary conditions
\begin{align}\label{eq:physicalaux}
u(x,t=0) = u_0; \qquad u(x=0,t) = u_\text{b} + s(t), \quad s(t) = a \sin \left( 2 \pi \nu t + \phi\right)
\end{align}
\end{subequations}
This problem describes, e.g., advection of a solute that undergoes linear decay; in this example $u$ represents {normalized} solute concentration, $v$ is {normalized} flow velocity along a streamline, and $k$ is the {normalized} reaction rate. In the simulations reported below, we set $v = 1$, $a=0.1$, $\nu=1$ and $\phi = 3 \pi /2$. 
In the first test, $k$ is a deterministic constant, while the uniform initial state $u_0$ and baseline state $u_\text{b}$ are random variables. In the other two tests, both $u_0$ and $u_\text{b}$ are deterministic, and $k$ is alternatively treated either as a random constant or as a spatially varying random field. 

In all three experiments, data sets $\mathbf d = \{ d_1, \cdots, d_{N_\text{meas}} \}$ are generated in accordance with~\eqref{eq:data} by adding Gaussian white noise, $\mathcal N(0, \sigma_\varepsilon)$, to a solution of~\eqref{eq:case0} with a given choice of model parameters. The likelihood function, $f_L(d_m | u(x,t)_m)$ with $ m =1, \cdots, d_{N_\text{meas}}$, is assumed to be Gaussian. 

The CDF equation for~\eqref{eq:case0} was derived, and the accuracy and robustness of the underlying closure approximations analyzed, in~\cite{boso-2014-cumulative} for the three scenarios described above. Appendix~\ref{app:MD} contains a brief summary of these results.


\subsection{Uncertain Initial and Boundary Conditions}
\label{sec:randominputs}


Let $u_0$ and $u_\text{b}$ be random uncorrelated random variables with (prior) PDFs $f_{u_0}(U_0)$ and $f_{u_\text{b}}(U_\text{b})$. Then the random initial and boundary states $u(x,t=0)$ and $u(x=0,t)$ are characterized by respective CDFs $F_0(U;\boldsymbol \varphi_0)$ and $F_\text{b}(U;t,\boldsymbol \varphi_\text{b})$ with shape parameters $\boldsymbol \varphi_0$ and $\boldsymbol \varphi_\text{b}$. In the absence of other sources of uncertainty, CDF $F_u(U;x,t)$ of the random state $u(x,t)$ in~\eqref{eq:case0} satisfies \emph{exactly} a PDE
\begin{subequations}\label{eq:CDF1d_case0}
\begin{align}
     \frac{\partial F_u}{\partial t} + \frac{\partial F_u}{\partial x} - k U \frac{\partial F_u}{\partial U} = 0 
\end{align}
subject to initial and boundary conditions
\begin{align}
  F_u(U;x,0) = F_0, \quad
  F_u(U;0,t) = F_\text{b}, \quad
  F_u(U_\text{min};x,t) = 0, \quad 
  F_u(U_\text{max};x,t) = 1.
\end{align}
\end{subequations}
This boundary-value problem admits an analytical solution, with either $F_0$ or $F_\text{b}$ that are propagated along deterministic characteristic lines. The dynamic manifold $\mathcal F$ of the resulting CDFs $F_u$ has coordinates $\tilde{\boldsymbol \varphi} = \{ x, t,\boldsymbol \varphi_0, \boldsymbol \varphi_\text{b}\}$. The analysis step of DA-MD takes place on this statistical manifold. Each measurement contributes to uncertainty reduction of either $\boldsymbol \varphi_0$ or $\boldsymbol \varphi_\text{b}$ (i.e., sharpens either $f_{u_0}$ or $f_{u_\text{b}}$), depending on the data location $(x,t)_m$. 
%
Half of these $N_\text{meas}$ measurements are collected at locations informing the initial condition, i.e., $(x/t)_m>1$), and the other half at locations informing the boundary condition, i.e., $(x/t)_m<1$. 


To verify the accuracy of DA-MD, we compare its predictions of the optimal parameters $\boldsymbol \varphi^{(N_\text{meas})}$ with those given by the  Bayesian posterior joint PDF 
\begin{align}\label{eq:bayesianinputs}
    \hat f_{u_0, u_\text{b}}(U_0, U_\text{b} | \mathbf d_{1:N_\text{meas}}) & = \hat f_{u_0}(U_0 | \mathbf d_{1:N_\text{meas}}) \hat f_{u_\text{b}}(U_\text{b} | \mathbf d_{1:N_\text{meas}}) \notag \\ 
    & \propto f_L(\mathbf d_{1:N_\text{meas}} | \mathbf u[(x,t)_{1:N_\text{meas}}; U_0, U_\text{b}] f_{u_0}(U_0) f_{u_\text{b}}(U_\text{b})\notag \\
    & \approx \prod_{m=1}^{N_\text{meas}} f_L(d_m | u[(x,t)_m; U_0, U_\text{b}] f_{u_0}(U_0) f_{u_\text{b}}(U_\text{b}).
\end{align}
%
To facilitate the Bayesian update, we take $F_{u_0}$ and $F_{u_\text{b}}$ to be Gaussian, fully specified by their respective means and standard deviations, $\boldsymbol \varphi_0 = \{ \mu_0, \sigma_0 \}$ and $\boldsymbol \varphi_\text{b} = \{ \mu_\text{b}, \sigma_\text{b} \}$. 
Then,~\eqref{eq:bayesianinputs} yields analytically-computable Gaussian posteriors $\hat f_{u_0}(U_0 | \mathbf d_{1:N_\text{meas}})$ and $\hat f_{u_\text{b}}(U_\text{b} | \mathbf d_{1:N_\text{meas}})$. In what follows, we compare those with the posterior parameters obtained via DA-MD, $\boldsymbol \varphi_0^{(N_\text{meas})}$ and  $\boldsymbol \varphi_\text{b}^{(N_\text{meas})}$, respectively. 
%
These posterior DA-MD parameters uniquely define the coefficients of the CDF equation~\eqref{eq:CDF1d_case0}, which then serves as an updated predictive tool. Equation~\eqref{eq:CDF1d_case0} has an analytical solution $F_u$ although, in general, numerical minimization in~\eqref{eq:loss2} needs to be employed to compute its approximation $F_\text{NN}$. 

\begin{figure}[htbp]
\centering
\includegraphics[width=\textwidth]{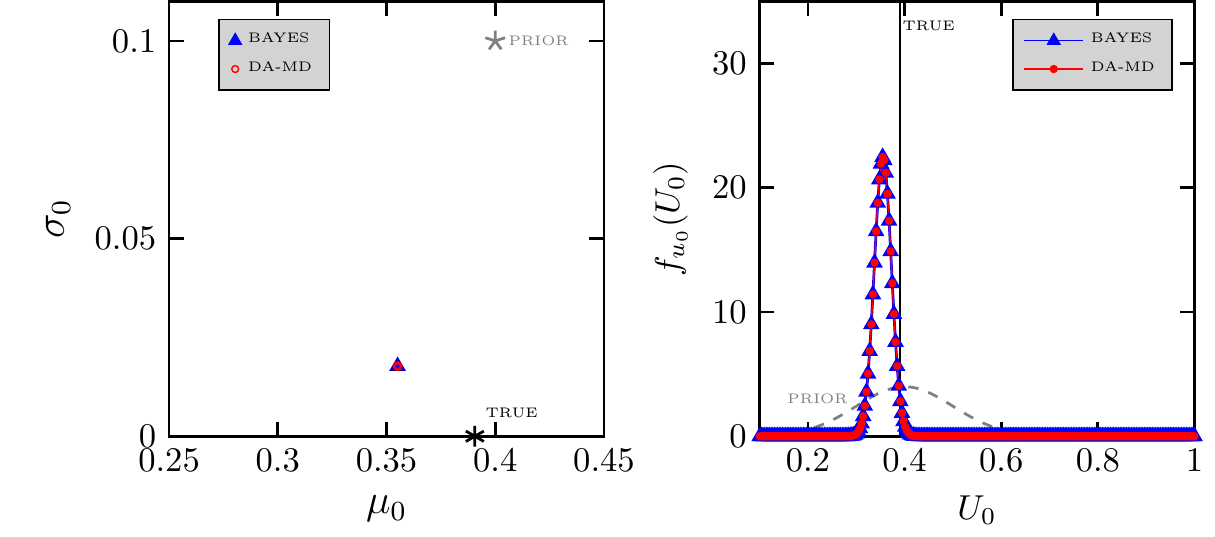}
\caption{Prior and posterior distributions for the initial state $u_0$ on the statistical manifold defined by the coordinates $\{ \mu_0, \sigma_0 \}$ representing the mean and standard deviation of a Gaussian distribution (left), and in the value space (right). 
The black asterisk in the left panel and the black vertical line in the right panel represent the true value ($u_0^{\text{true}} = 0.391$)
, for which a Gaussian PDF degenerates into the  Dirac distribution (delta function). The grey star (left) and the gray dashed line (right) represent a prior distribution ($\mu_0^{\text{prior}} = 0.4$, $\sigma_0^{\text{prior}} = 0.1$). 
The blue triangle (left) and line (right) identify the Bayesian solution, whereas the corresponding red symbols and lines  identify the DA-MD solution. Parameters are set to $k = 1$, $\sigma_\varepsilon = 0.04$ and $N_\text{meas}=20$.
}
\label{fig:posterior_case0}
\end{figure}

Figure~\ref{fig:posterior_case0} exhibits the prior and posterior distributions for $u_0$ (those for $u_b$ behave similarly) computed with the alternative data assimilation strategies. The left panel represents these distributions as coordinates ($\mu_0, \sigma_0$) on the statistical manifold of Gaussian distributions, whereas the right panel shows them as PDFs in the value space $\Omega_{U_0}$. The Bayesian update and the DA-MD approach 
yield almost identical results after assimilation of the same set of measurements, sharpening the distribution of the parameters around the true value. 

Similar to the left panel in \cref{fig:posterior_case0}, the prior and posterior CDFs of the state variable $u(x,t)$, both obeying the CDF equation~\eqref{eq:CDF1d_case0}, are represented as points on the statistical manifold $\mathcal F$ with coordinates $(x,t,\boldsymbol \varphi^{(0)})$ and $(x,t,\boldsymbol \varphi^{(N_\text{meas})})$, respectively.
The 
amount of information used during the analysis and transferred from the measurements to the conditional predictions can be thought of as the distance between these points: the information gain from prior to posterior is quantified by the KL divergence between these distributions (\cref{sec:analysis}\ref{sec:ig}). 
For the same prior and the same observations, DA-MD and the Bayesian update yield almost identical KL discrepancies. Moreover, $D_{\text{KL}}$ does not vary within the assimilation regions, i.e., it remains constant in the regions of the space-time domain where $F_u$ depends on either $\boldsymbol \varphi_0$ or $\boldsymbol \varphi_b$. The KL divergence also allows one to compare the informational gain from different sets of observations: doubling the number of measurements from $N_{\text{meas}} = 20$ to $N_{\text{meas}} = 40$ yields, in the assimilation regions informed by either the initial or the boundary conditions, a gain in KL terms of 7\% and 9\%, respectively.

\subsection{Uncertain Reaction Rate}
\label{sec:randomk}

In the following two test-cases, we treat the uncertain coefficient $k$ in~\eqref{eq:case0} first as a random constant and then as a random field. The auxiliary variables $u_0$ and $u_{\text b}$ in~\eqref{eq:physicalaux} are taken to be deterministic, so that the CDF equation~\eqref{eq:CDF} is subject to initial and boundary conditions
\begin{align*}
    F_u(U;x,0) = \mathcal H(U - u_0), \qquad F_u(U;b,t) = \mathcal H(U - u_{\text b} - s(t)).
\end{align*}

\subsubsection{Random variable}

The coefficients~\eqref{eq:exactclosures_deterministicdynamics} in the CDF equation~\eqref{eq:CDF} take the form (\cref{app:MD})
\begin{align}\label{eq:closurek}
\boldsymbol{\mathcal Q} = \begin{pmatrix}
1 \\
- \langle k \rangle U - \dfrac{\sigma_k^2 U}{ \langle k \rangle } [1- \text{e}^{\langle k \rangle t^*} ]
\end{pmatrix}, \qquad
  \boldsymbol{\mathcal D} =  \begin{pmatrix}
0 & 0 \\
0 & - \dfrac{\sigma_k^2 U^2}{\langle k \rangle } [1- \text{e}^{\langle k \rangle t^*} ]
\end{pmatrix},
\end{align}
where $t^* (U,x,t)= \min \{ t,x,\langle k \rangle^{-1} \ln ( U_{\text{max}}/U )\}${, and $\langle k \rangle$ and $\sigma_k$ are the ensemble mean and standard deviation of $k$, respectively}.
The coordinates of the dynamic manifold $\mathcal F$ of the approximated CDFs $F_u$ are $\tilde{\boldsymbol \varphi} = \{ x,t,\langle k \rangle, \sigma_k\}$. 
The CDF equation is solved via finite volumes (FV) using the Fipy solver \cite{guyer2009fipy}, setting the discretization elements to $\Delta t = 0.01$, 
$\Delta x = L/200$ and $\Delta U = ( U_{\text{max}}-U_{\text{min}} )/128$, with domain size defined by $L=1$, $U_\text{min} = 0$ and $U_\text{max} = 1$.

Minimization of \eqref{eq:loss2} is done using the L-BFGS-B method implemented in \texttt{TensorFlow}~\cite{tensorflow2015-whitepaper} {with a convergence threshold for the loss function value of $10^{-3}$}. The solution of the CDF equation~\eqref{eq:CDF}, whose coefficients are given by~\eqref{eq:closurek}, is represented by a fully connected NN with fixed architecture (9 layers, 20 nodes per hidden layer) and a sigmoidal activation function (hyperbolic tangent). Weights and biases of the NN are initialized at the beginning of the sequential procedure by approximating a solution of the CDF equation with prior statistical parameters $\boldsymbol \varphi^{\text{(0)}}$. Successive iterations are initialized with weights and biases from the previous step. This procedure considerably accelerates the identification of the target parameters. Zero residuals are enforced at $N_{\text{res}} = 792$ 
locations within the space-time domain, whereas initial and boundary conditions are imposed at $N_{\text{aux}} = 406$ locations. Furthermore, we enforce non-negativity of $\sigma_k$. 

\begin{remark}
The FV approximation is used to construct the observational CDFs, whereas the NN approximation is used on sparse set of points for numerical gradient-based minimization. The NN surrogate solution of the CDF equation~\eqref{eq:CDF} could also be used as a prior for the next assimilation step, with the advantage of being virtually free of artificial diffusion and with no theoretical limitation on the number of dimensions. This is not exploited further in this work, as research on the use of physics-informed NN to solve PDEs is not yet mature. Nevertheless, it has been shown to yield accurate identification of PDE parameters~\cite{raissi2019physics} and to reproduce qualitatively actual PDE solutions.
\end{remark}

\begin{figure}[htbp]
\centering\includegraphics[width=\textwidth]{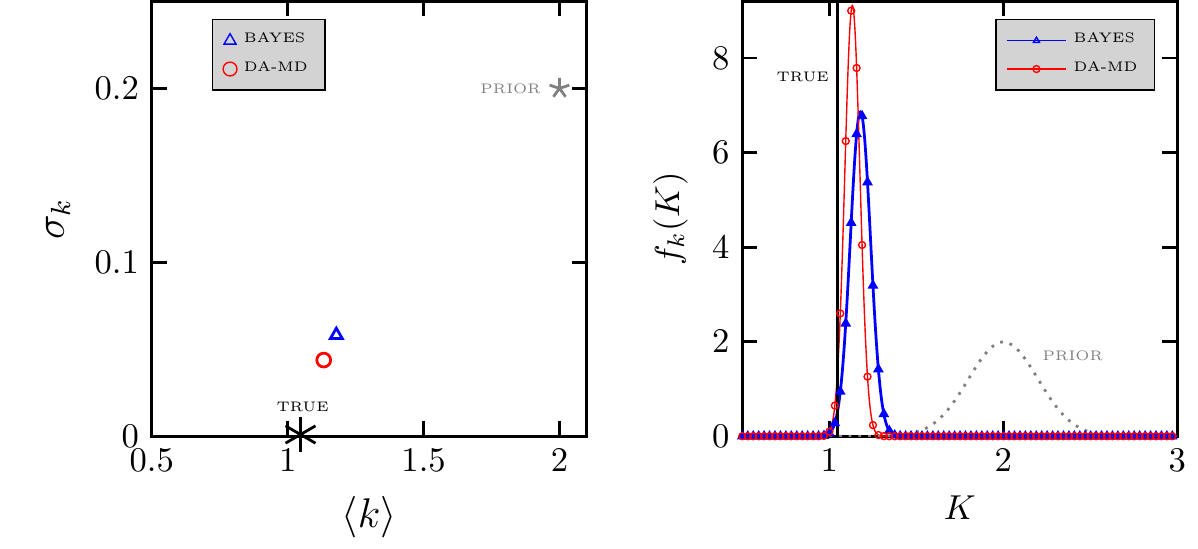}
\caption{Prior and posterior PDFs of the random variable $k$ shown (left) on the statistical manifold defined by coordinates $\{ \langle k \rangle,\sigma_k \}$ representing the mean and standard deviation of $k$, and (right) in the value space $\Omega_K$. The black asterisk in the left panel and the black vertical line in the right panel represent the true value $(k^{\text{(TRUE)}}=1.047)$, for which a Gaussian PDF degenerates into the Dirac distribution (delta function). The gray star (left) and the gray dashed line (right) represent a prior distribution $(\langle k \rangle^{\text{prior}} = 2, \sigma_k^{\text{prior}} = 0.2)$. The blue triangle (left) and line (right) identify the Bayesian solution, whereas the corresponding red circles and lines identify the DA-MD solution. Parameters are set to $u_0=0.4$, $u_{\text b}=0.5$, $\sigma_\varepsilon = .02$, $N_{\text{meas}}=20$.
}
\label{fig:posterior_case1}
\end{figure}

We compare the DA-MD estimate of the PDF of the model parameter $k$ with the Bayesian posterior PDF of $k$. The latter is obtained analytically by assuming a Gaussian prior $f_k(K)$ and taking advantage of the analytical solution of \eqref{eq:case0}
\begin{align}
\hat f_k(K|\mathbf d_{1:N_{\text{meas}}}) & \propto f_L(\mathbf d_{1:N_{\text{meas}}} | \mathbf u[(x,t)_{1:N_{\text{meas}}}; K ]) f_k(K) \nonumber\\
& \approx \sum_{m=1}^{N_{\text{meas}}} f_L (d_m | u \left[ (x,t)_m; K\right]) f_k(K).
\end{align}
The Bayesian and DA-MD posterior (and prior) PDFs of the random reaction rate $k$ are presented in \cref{fig:posterior_case1}. The right panel shows these densities in the value space $\Omega_K$ of $f_k(K)$, whereas the left panel represents the state distributions as points on the dynamic manifold $\mathcal F$.  The Bayesian update is optimal and analytical. Its sole source of error stems from the calculation of the normalization constant via numerical integration; as such it is treated as a benchmark in this comparison. On the contrary, DA-MD is based on a series of approximations (closures for the CDF equation, FV and NN solutions of the CDF equation, numerical minimization of the loss function). Nevertheless, DA-MD yields an updated posterior which is close to the Bayesian one. The DA-MD posterior is sharper than the Bayesian posterior; this might be due to the effect of numerical diffusion that artificially smears the CDF profiles computed as a solution of the CDF equation.  
Convergence of the DA-MD is slow, but its computational time is not expected to scale with the dimensionality of the problem (e.g., when dealing with random parameter fields). This flexibility represents a major advantage of the proposed procedure versus Bayesian inference, and it is explored in a more challenging scenario in the following section. 

\begin{figure}[htbp]
\centering\includegraphics[width=\textwidth]{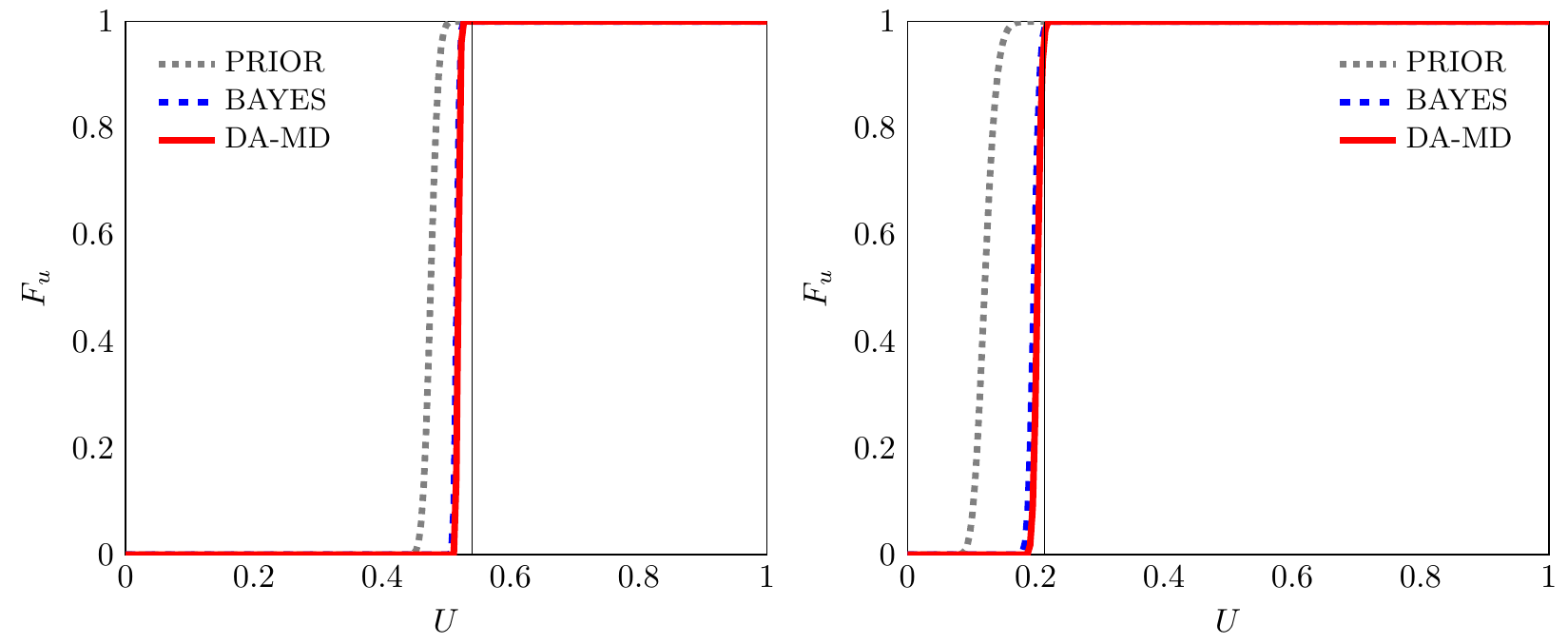}
\caption{Prior and posterior CDF profiles $F_u(U)$ for the random $k$ case at the final assimilation time $t_M$ at two spatial locations: $x=0.1$ (left) and $x=0.8$ (right). The vertical black line represents the true solution, for which the CDF degenerates into a Heaviside step function. The dotted grey line represents the prior distribution with parameters $(\langle k \rangle, \sigma_k) = \{ 2,0.2\}$; the dashed blue line and the solid red line represent the posterior distribution computed with updated Bayesian and DAMD parameters, respectively. The updated parameters $\langle k \rangle$ and $\sigma_k$ are those represented in \cref{fig:posterior_case1}.  The remaining parameters are set to $k^{\text{(true)}} = 1.047$,$u_0=0.4$, $u_{\text b}=0.5$, $\sigma_\varepsilon = 0.02$, $N_{\text{meas}}=20$, $t_M=0.6$.
}
\label{fig:CDFs_case1}
\end{figure}

The prior and posterior CDFs of $u$, $F_u(U;\cdot)$ at the final assimilation time $t_M$ are shown in \cref{fig:CDFs_case1}. The posterior CDFs, for both the Bayesian and DA-MD assimilation, provide a state prediction that is closer than the prior CDF to its true value thanks to a more accurate parameter identification (shown in \cref{fig:posterior_case1}). The value of measurements is evaluated in terms of their impact on the shape of the CDF at the measurement locations, and quantified by the KL divergence from the posterior to the prior. In this example, all locations exhibit the same information gain quantified by the KL divergence going from the posterior to the prior. That is because of the analytical one-to-one relation between $k$ and $u(\mathbf x,t)$.


\subsubsection{Random field}

Keeping all other conditions and settings unchanged, we now consider a spatially varying uncertain parameter $k(x)$. It is treated as a second-order stationary (statistically homogeneous) multivariate Gaussian random field with constant mean $\langle k \rangle^{\text{(true)}}$ and standard deviation $\sigma_k^{\text{(true)}}$; its two-point autocovariance function $C_k^{\text{(true)}} ( x -  x') $ has either zero correlation length (i.e., uncorrelated random field or white noise),
\begin{align*}
C_k^{\text{(true)}} ( x -  x') = \sigma_k^2 \delta{\left(  x -  x' \right)},
\end{align*}
or a finite correlation length $\lambda_k^{\text{(true)}}$,
\begin{align*}
C_k^{\text{(true)}} ( x -  x') = \sigma_k^2 \exp{\left( -| x -  x'| / \lambda_k^{\text{(true)}} \right)}.
\end{align*}
One realization with the chosen statistical parameters represents the reference random field $k^{\text{(true)}}(x)$, which was used to construct synthetic observations via the FV solution of~\eqref{eq:conslaw} with~\eqref{eq:case0}. 

The coefficients~\eqref{eq:exactclosures_deterministicdynamics} in the CDF equation~\eqref{eq:CDF} now take the form (\cref{app:MD})\begin{align}\label{eq:closureterms_uncorr}
\boldsymbol{\mathcal Q} =  
\begin{pmatrix} 1 \\ 
- \langle k \rangle U - \dfrac{\sigma_k^2 U}{2} 
\end{pmatrix}, \qquad
 \boldsymbol{\mathcal D} = 
 \begin{pmatrix}
 0 & 0 \\ 
 0 & \dfrac{\sigma_k^2 U^2 }{2}
 \end{pmatrix},
\end{align}
if $k(x)$ is white noise, or 
\begin{align}\label{eq:closureterms_finite}
 \boldsymbol{\mathcal Q} = 
 \begin{pmatrix} 1 \\ - \langle k \rangle U - \dfrac{\sigma_k^2 U}{\alpha} [\text{e}^{ \alpha t^*} -1 ] 
 \end{pmatrix}, \qquad
 \boldsymbol{\mathcal D} = 
 \begin{pmatrix}
  0 & 0 \\ 
  0 & \dfrac{\sigma_k^2 U^2}{\alpha} [ \text{e}^{ \alpha t^* } -1] 
\end{pmatrix}
\end{align}
with $\alpha = \langle k \rangle - 1 / \lambda_k$ and $t^*(U,x,t) = \min \{ t,x,\langle k \rangle^{-1} \ln (U_{\text{max}} / U) \}$, if $k(x)$ has the exponential correlation $C_k$. 
The corresponding dynamic manifolds have either the coordinates $\tilde{\boldsymbol \varphi} = \{ x,t,\langle k \rangle, \sigma_k\}$ or the coordinates $\tilde{\boldsymbol \varphi} = \{ x,t,\langle k \rangle, \sigma_k, \lambda_k\}$, respectively. 

Unlike Bayesian update, which identifies the $k$ values at each spatial location with a consequent dramatic increase of the dimensionality of the target joint posterior PDF, DA-MD focuses on a finite set of parameters $\boldsymbol {\varphi}$ (the mean $\langle k \rangle$, the standard deviation $\sigma_k$ and, in the correlated case, the correlation length $\lambda_k$); its computational cost is comparable to that for the constant random parameter case. We compare the updated DA-MD parameters with an approximation of the Bayesian posterior, since the number of random parameters and the nonlinearity of the problem prevent analytical treatment. 

We employ a standard ensemble Kalman filter (EnKF) \cite{wikle2007bayesian, evensen2009data} for the update of the random field $k(x)$, discretized into $N_x$ point values, with ensemble size $N_{\text{ens}}$. EnKF requires multiple solutions of the physical model, which typically require special numerical treatment because of the high spatiotemporal  variability of the model parameters. The choice of a spatial resolution poses another difficulty because the correlation length of the target random field is a-priori unknown. This increases the numerical complexity of EnKF, to the advantage of MD. To focus on the data assimilation aspect of the problem, we solve both the physical model and the CDF equation on the same grid and with the same FV numerical solver, thus taking advantage of MD's lower numerical complexity. Update is done sequentially for DA-MD, and recursively for EnKF \cite[and references therein]{crestani2013ensemble}, i.e., at each assimilation step the ensemble members are forecast from the initial time to the current assimilation time. In both cases (EnKF and DA-MD), we assimilate the same $N_{\text{meas}}$ measurements collected at two spatial locations, $x=0.1$ and $x=0.8$, in ten separate temporal instances, $t = \{ 0.15, 0.2, \dots, 0.6 \}$. 

\begin{figure}[htbp]
\centering\includegraphics[width=\textwidth]{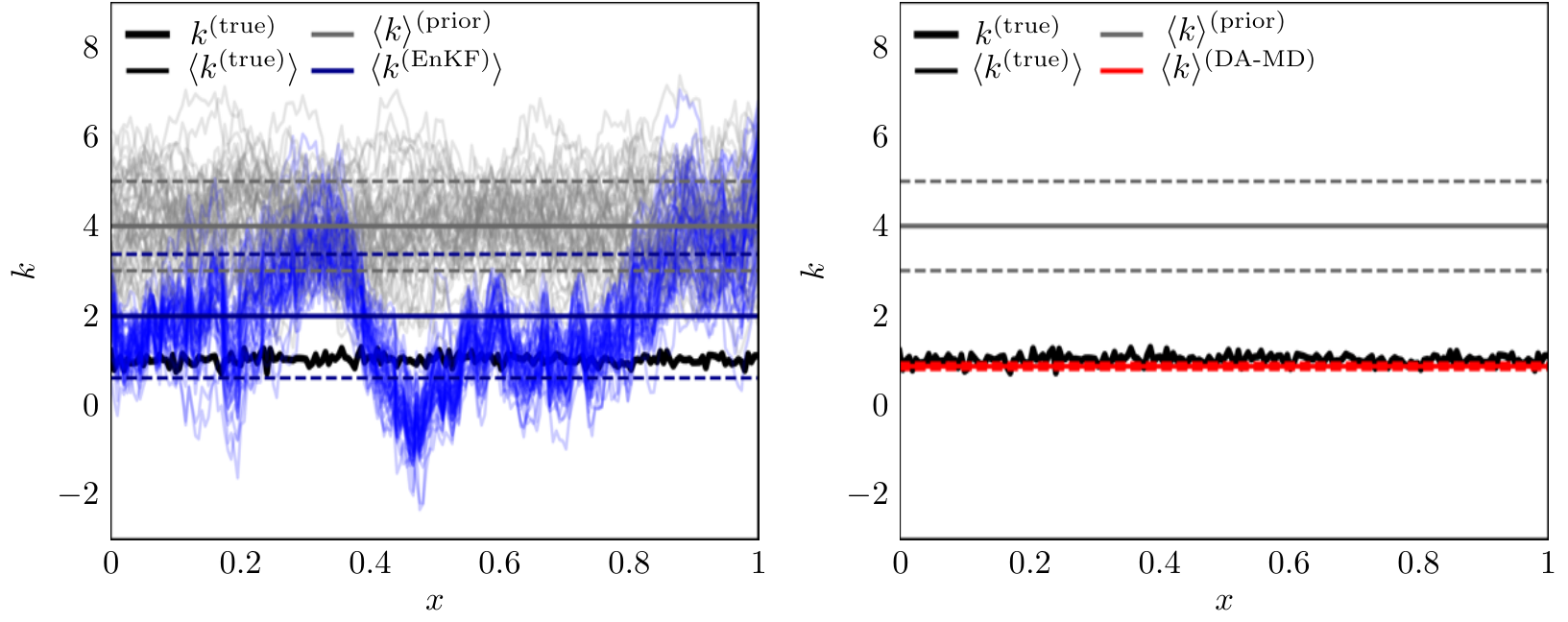}
\caption{
Parameter identification for the uncorrelated $k(x)$ field. Both panels contain the true field, $k^{\text{(true)}}(x)$, in black, and the prior field moments (grey lines). Both the prior and posterior random fields are defined by their mean value $\langle k \rangle$ (solid line), and a buffer region with half-width equal to the standard deviation (dashed lines). For the EnKF (left), both the prior and posterior ensemble members are represented. 
Posterior values are $\langle k \rangle^{\text{(DA-MD)}} = 0.86$, $\sigma_k^{\text{(DA-MD)}} = 0.07$, 
$\langle \overline k^{\text{(EnKF)}} \rangle = 1.7$, $\overline \sigma_k^{\text{(EnKF)}} = 1.19$, $\overline k^{\text{(true)}} = 1.01$, 
$\overline \sigma_k^{\text{(true)}} = 0.1$, where $\overline \cdot$ represents the spatial average. Parameters are set to $u_0=0.4$, $u_{\text{ b}} = 0.5$, $N_{\text{ens}} = 50$, $N_{\text{meas}}=20$, $\sigma_\varepsilon = 0.02$, $N_x = 200$, $\boldsymbol \varphi^{(0)} = \{ \langle k \rangle, \sigma_k \}^{\text{(prior)}} = (4,1)$.}
\label{fig:posterior_statpars_uncorr}
\end{figure}

\begin{figure}[htbp]
\centering\includegraphics[width=\textwidth]{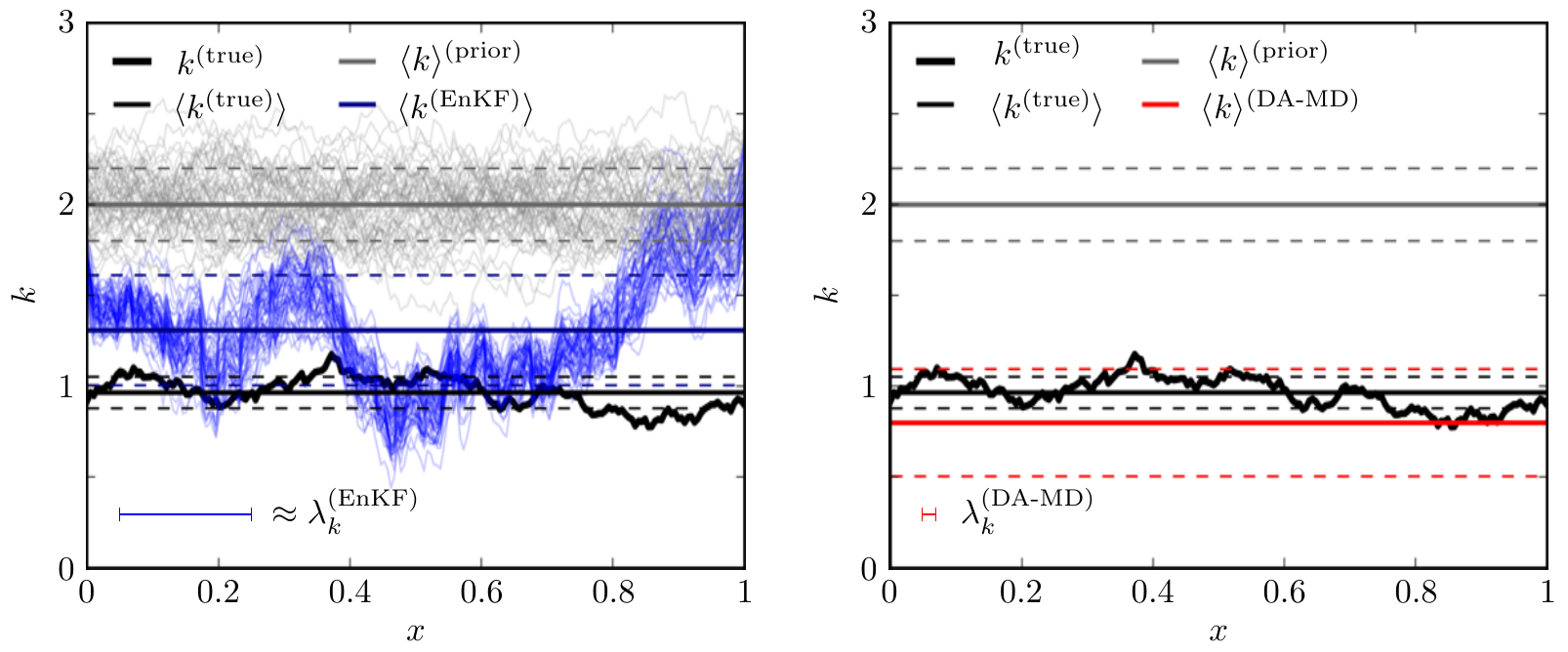}
\caption{{Parameter identification for the correlated field $k(x)$ via EnKF (left) and DA-MD (right). Both panels contain the true field, $k^{\text{(TRUE)}}(x)$, in black, and the prior field (grey lines). Both the prior and posterior fields are defined by their mean value $\langle k \rangle$ (solid line), and a buffer region with half-width equal to the standard deviation value (dashed lines). An estimate of the posterior correlation length is in the bottom left corner of both panels. For EnKF (left panel), both the prior and posterior ensemble members are also represented. Posterior values are $\langle k \rangle^{\text{(DA-MD)}} = 0.80$, $\sigma_k^{\text{(DA-MD)}} = 0.30$, $\lambda_k^{\text{(DA-MD)}} = 0.013$, $\langle \bar k^{\text{(EnKF)}} \rangle = 1.23$, $\bar \sigma_k^{\text{(EnKF)}} = 0.33$, $\bar k^{\text{(true)}} = 0.96$, $\bar \sigma_k^{\text{(true)}} = 0.09$, $\lambda_k^{\text{true}} = 0.3$, where $\overline \cdot$ represents the spatial average. Parameters are set to $u_0=0.4$, $u_{\text b} = 0.5$, $N_{\text{ens}} = 50$, $N_{\text{meas}}=20$, $\sigma_\varepsilon = 0.01$, $N_x = 200$, $\boldsymbol \varphi^{(0)} = \{ \langle k \rangle, \sigma_k, \lambda_k \} ^{\text{(prior)}} = \{ 2,0.2,0.2 \}$.}}
\label{fig:posterior_statpars_corr}
\end{figure}

Figures~\ref{fig:posterior_statpars_uncorr} and~\ref{fig:posterior_statpars_corr} exhibit the EnKF and DA-MD posterior random fields  for the uncorrelated and correlated cases, respectively. When the true field $k(x)^{\text{(true)}}$ is white noise, DA-MD accurately identifies the updated mean $\langle k \rangle^{(\text{DA-MD})}$, but underestimates the value of $\sigma_k^{\text{(DA-MD)}}$. The latter might be due to the impact of artificial diffusion on the solution of the CDF equation used as a prior in the DA-MD procedure. EnKF yields a wider posterior estimate for $k$, with spatial averages for the mean $\langle k \rangle^{\text{(EnKF)}}$ and the standard deviation $\sigma_k^{\text{(EnKF)}}$ that are further away from the spatial averages of the moments of the true field (values in the caption). 

\begin{figure}[htbp]
\centering\includegraphics[width=\textwidth]{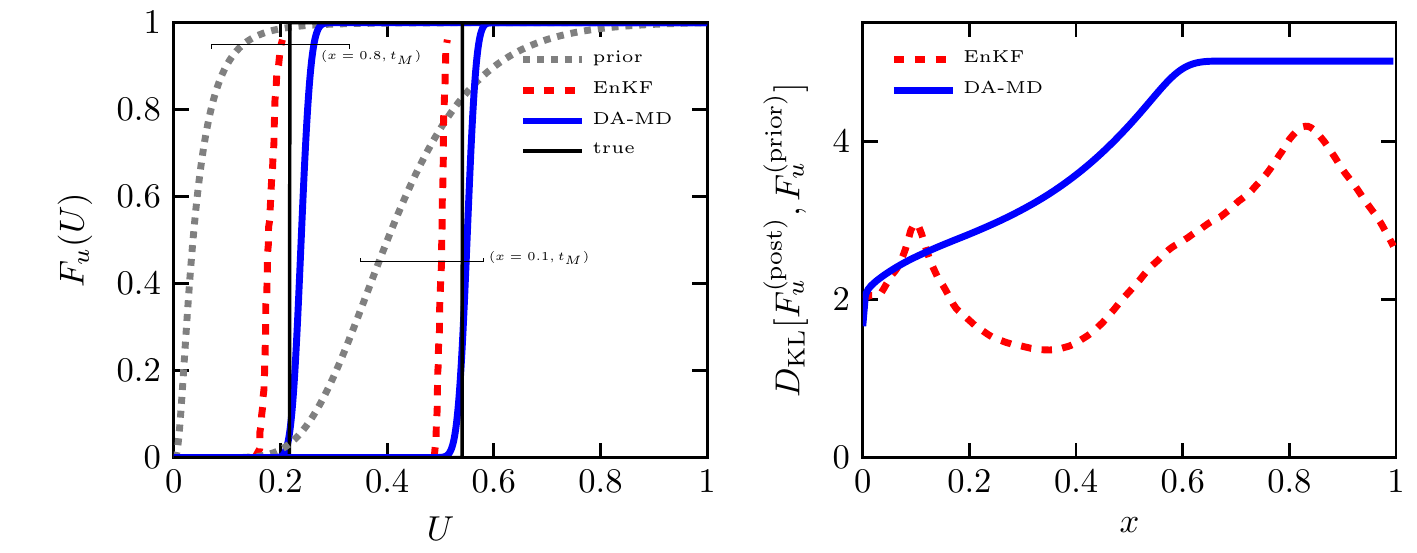}
\caption{{Prior and posterior CDFs (left) and corresponding KL divergence $D_{\text{KL}}$ (right) for the uncorrelated $k(x)$ field. The CDF profiles (left) are computed at $(x,t)=(0.1,t_M)$ and $(x,t)=(0.8,t_M)$ as a solution of the CDF equation with prior $\boldsymbol \varphi^{(0)}$ and posterior $\boldsymbol \varphi^{(N_{\text{meas}})}$ parameters (dotted grey and solid blue lines, respectively). The CDFs from EnKF (dashed red line) are computed as an empirical distribution of the ensemble members. The true solution is plotted as a Heaviside function centered on the true value $u^{\text{(true)}}(x,t)$ (black thin line), $\mathcal H(U - u^{\text{(true)}}(x,t))$. The selected coordinates for the profiles ($x=0.1$ and $x=0.8$) correspond to measurement locations. For both DA-MD and EnKF, the KL divergence $D_{\text{KL}}$ between the posterior distribution and the prior distribution is computed as a function of $x$ at time $t_M$. 
Parameters are set to $u_0=0.4$, $u_{\text b} = 0.5$, $N_{\text{ens}} = 50$, $N_{\text{meas}}=20$, $\sigma_\varepsilon = 0.02$, $N_x = 200$, $\Delta x = 1.6 \cdot 10^{-3}$, $\Delta U = 8.3 \cdot 10^{-4}$, $\Delta t = 10 ^{-3}$, $\boldsymbol \varphi^{(0)} = \{ \langle k \rangle, \sigma_k \} ^{\text{(prior)}} = \{ 4,1\}$, $t_M = 0.6$.}}
\label{fig:posterior_profiles_uncorr}
\end{figure}

Despite these differences in reconstruction of the statistical properties of the posterior $k(x)$, both assimilation techniques yield a posterior prediction of $F_u(U)$ that approaches the true value of the solution (the left panel in  \cref{fig:posterior_profiles_uncorr}). The information gain from the measurements is quantified in terms of the KL divergence for both DA-MD and EnKF (the right panel in  \cref{fig:posterior_profiles_uncorr}) at time $t_M$. MD densities (both the prior and the posterior) are calculated via finite differences from the solution of the CDF equation, whereas EnKF densities are computed via Kernel Density Estimation with Gaussian kernel and Scott's bandwidth, using the ensemble members as data points.  
Our results suggest DA-MD extracts more information than EnKF from the same set of measurements in the current configuration at almost all values of $x$, as is also reflected in an accurate characterization of  the posterior $k$ field. Observations collected at $x>t_M$ (the region where characteristic lines originate from the initial conditions) are more informative for DA-MD assimilation. The KL divergence for EnKF highlights the locations of more informative measurements, displaying two distinctive peaks.

Figure~\ref{fig:posterior_statpars_corr} exhibits the results of a similar analysis for the correlated field $k(x)^{\text{(true)}}$. DA-MD posterior estimates of the mean and standard deviation of $k$ are closer to the averaged statistical properties of the true field than EnKF estimates  are (values are in the figure caption). DA-MD underestimates the spatial correlation length $\lambda_k$, whereas the identification of $\lambda_k$ via EnKF is inconclusive as the semivariogram for $k(x)$ does not develop a sill. We identify an intermediate plateaux and assume the corresponding lag value as the updated correlation length for the field. The semivariogram is computed using the posterior ensemble member values, and is shown in \cref{fig:posterior_statpars_corr}(right). The corresponding state CDFs $F_u(U;x,t)$ are plotted in \cref{fig:posterior_statpars_corr}(left) in two representative sections that correspond to measurement locations. Both DA-MD and EnKF yield a posterior state CDF $F_u$ considerably closer to the true value than the prior distribution.

\begin{figure}[htbp]
\centering\includegraphics[width=\textwidth]{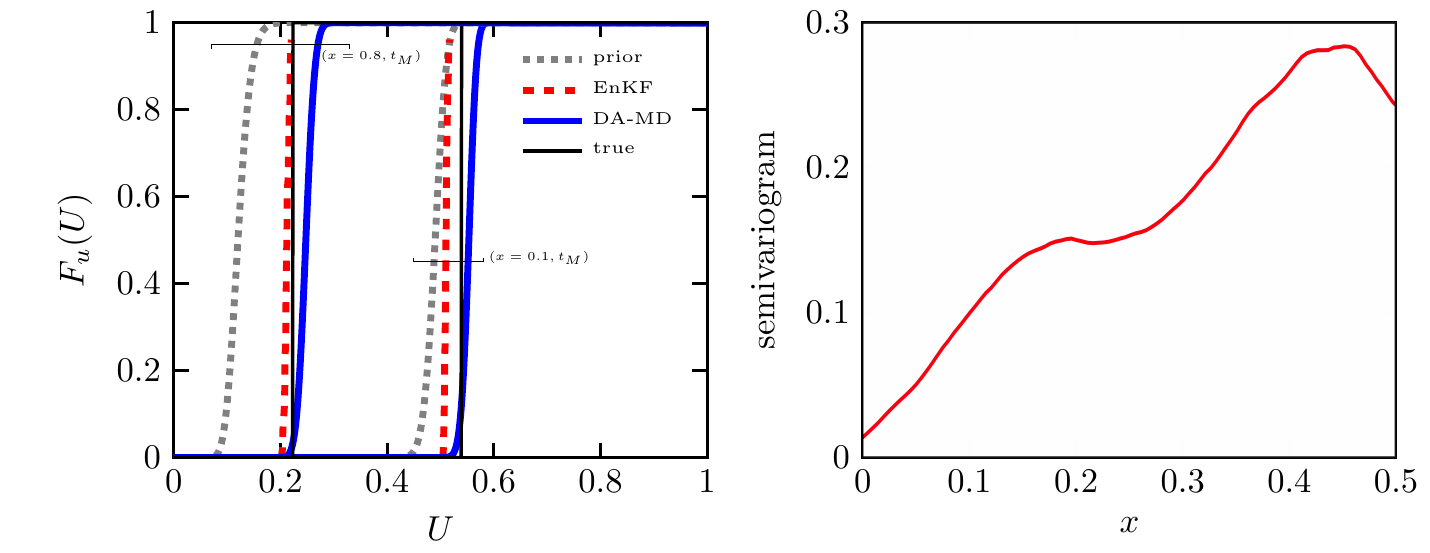}
\caption{{\textit{Left}: Prior and posterior CDFs of the correlated field $k(s)$. The CDF profiles are computed at $(x,t)=(0.1,t_M)$ and at $(x,t) = (0.8,t_M)$ as a solution of the CDF equation with prior $\boldsymbol \varphi^{(0)}$ and posterior $\boldsymbol \varphi^{(N_{\text{meas}})}$ parameters (dotted grey and solid blue lines, respectively). The CDFs from EnKF (dashed red line) are computed as an empirical distribution of the ensemble members. The true solution is plotted as a Heaviside function centered on the true value $u^{\text{(true)}}(x,t)$ (black thin line), $\mathcal H(U - u^{\text{(true)}}(x,t))$. \textit{Right:} Semivariogram for the EnKF posterior ensemble members. Parameters are set to $u_0=0.4$, $u_{ \text b} = 0.5$, $N_{\text{ens}} = 50$, $N_{\text{meas}}=20$, $\sigma_\varepsilon = 0.01$, $N_x = 200$, $\Delta U = 3.75 \cdot 10^{-4}$, $\Delta x = 1.6 \cdot 10^{-3}$, $\Delta U = 8.3 \cdot 10^{-4}$, $\Delta t = 10 ^{-3}$, $\boldsymbol \varphi^{(0)} = \{ \langle k \rangle, \sigma_k,\lambda_k \} ^{\text{(prior)}} = (2,0.2,0.2)$, $t_M = 0.6$.}}
\label{fig:posterior_profiles_corr}
\end{figure}


\section{Summary and Future Work}
\label{sec:conclusion}

We proposed a novel methodology for parameter estimation that leverages the method of distributions (MD) for both the forecast and analysis steps. Reduction of uncertainty on model parameters is recast into a problem of identification of closure parameters for the CDF equation, expressing the space-time evolution of uncertainty for the model output. Specifically,  we identify the parameters in the CDF equation \eqref{eq:CDF} which yield an estimate in the measurement locations as close as possible to the state distribution. This is expressed by an observational Bayesian posterior in that specific location, which is obtained combining the data model and the physically-based prior. The procedure is done sequentially, progressively updating the parameters of the CDF equation as more measurements are assimilated. We demonstrated that our method reproduces Bayesian posteriors in scenarios where Bayesian inference can be performed analytically, and ameliorates parameter identification when compared to ensemble Kalman filter (as an approximation of Bayesian update) in cases where Bayesian inference is elusive.

This work opens multiple possible research venues. In particular, we plan to i) explore the construction of novel data-driven closure approximations for MD; ii) investigate the use of novel ML techniques for more efficient optimization and/or solution of PDEs; iii) introduce multi-point statistics.

\appendix

\section{Derivation of CDF Equations}
\label{app:MD}

MD commences by defining a so-called raw CDF $\pi(U;\mathbf x,t) \equiv \mathcal H(U - u(\mathbf x,t))$, where $\mathcal H(\cdot)$ is the Heaviside function. Let $f_u(U;\mathbf x,t)$ denote the single-point PDF of $u(\mathbf x,t)$. Then it follows from the definition of the ensemble mean $\mathbb E[\cdot] \equiv \langle \cdot \rangle$ that 
\begin{align} \label{eq:mean_pi}
   \mathbb E[ \pi(U;\mathbf x,t) ] &\; = \int_{U_\text{min}}^{U_\text{max}} \mathcal H(U - \mathcal U) f_u (\mathcal U; \mathbf x,t) \text d \mathcal U = \int_{U_\text{min}}^U \mathcal H(U - \mathcal U) f_u (\mathcal U; \mathbf x,t) \text d \mathcal U \nonumber\\
   &\; = F_u(U;\mathbf x,t).
\end{align}
Other useful properties of $\pi$ are
\begin{align}
    \frac{\partial \pi}{\partial t} = \frac{\partial \pi}{\partial u} \frac{\partial u}{\partial t} = - \frac{\partial \pi}{\partial U} \frac{\partial u}{\partial t} \quad \text{and} \quad \nabla \pi = - \frac{\partial \pi}{\partial U} \nabla u.
\end{align}
Accounting for these properties, multiplication of~\eqref{eq:conslaw} by $- \partial_U \pi$ yields
\begin{align}\label{eq:pi}
    \frac{\partial \pi}{\partial t} + \dot{\mathbf q}(U) \cdot \nabla \pi + r(U) \frac{\partial \pi}{\partial U} = 0,
\end{align}
{where $\dot{\mathbf q} = \text d \mathbf q(U) / \text d U$.}
This equation is exact as long as solutions of \eqref{eq:conslaw}, $u(\mathbf x,t)$, are smooth (do not develop shocks) for each realization of random parameters $\tilde{\boldsymbol\theta}$. It is subject to initial and boundary conditions derived from the initial and boundary conditions of the physical problem, and to $\pi(U = U_{\text{min}};\mathbf x,t) = 0$ and $\pi(U = U_{\text{max}};\mathbf x,t) = 1$. 

In the absence of uncertainty, \eqref{eq:pi}  is deterministic and equivalent to \eqref{eq:conslaw}; the model output $u(\mathbf x,t)$ can be recovered from $\Pi(U,\mathbf x,t)$ by integration. In the presence of uncertainty affecting the parameters and the auxiliary inputs, it follows from~\eqref{eq:mean_pi} that the ensemble average of~\eqref{eq:pi} is 
\begin{align}\label{eq:piaveraged}
    \frac{\partial F_u}{\partial t} + \langle \dot{\mathbf q} (U; \boldsymbol \theta_q) \cdot \nabla \pi \rangle + \langle r(U;\boldsymbol \theta_r) \frac{\partial \pi} {\partial U} \rangle.
\end{align}
If the model parameters $\boldsymbol \theta$ are deterministic, then so is the evolution dynamics and uncertainty in predictions of $u(\mathbf x,t)$ is solely due to uncertainty in the initial and the boundary conditions. In that case,~\eqref{eq:piaveraged} gives an exact CDF equation,
\begin{align}
    \frac{\partial F_u}{\partial t} +  \dot{\mathbf q} (U; \boldsymbol \theta_q) \cdot \nabla F_u  +  r(U;\boldsymbol \theta_r) \frac{\partial F_u}{\partial U}.
\end{align}
Otherwise, closure approximations are necessary to obtain a workable expression for the undefined terms in \eqref{eq:piaveraged}. These expressions depend on the closure strategy and on the functional form of $\mathbf q$ and $r$. 

To be specific, we set $\mathbf q(u) = \mathbf v(\mathbf x) u$ and $r(u) = k r_\alpha(u; \alpha, u_\text{eq}) = k \alpha \left( u_\text{eq}^\alpha - u^\alpha \right)$. Here $\mathbf v(\mathbf x)$ is the divergence-free velocity, $\nabla \cdot \mathbf v = 0$, of steady incompressible flow; and $\alpha \in \mathbb N^+$ is the order of an equilibrium reaction with reaction rate $k(\mathbf x)$, which drives the system towards its equilibrium state $u_{\text{eq}}$. An analogous system was studied in detail in \cite{boso-2014-cumulative, venturi2013exact}. In what follows we summarize the closure approximations developed in these works for the case of deterministic $\mathbf v(\mathbf x)$ and random $k(\mathbf x)$.

We use the Reynold decomposition to represent random quantities as the sum of their respective means and zero-mean fluctuations around these means,  
\begin{align}
    k = \langle k \rangle + k', \qquad \pi = F + \pi'. 
\end{align}
A first-order (in the variance $\sigma_k^2$ of stationary random fluctuations $k'$) approximation of~\eqref{eq:piaveraged} takes the form of \eqref{eq:CDF} with the coefficients~\cite{boso-2014-cumulative} 
\begin{align}\label{eq:closuresG}
    & \mathcal Q_i =  v_i (\mathbf x) , \quad i=1,\dots, d \notag \\
    & \mathcal Q_{d+1} \approx \langle k \rangle r_\alpha (U) + \int_0^t \int_{\tilde \Omega} G(\mathbf x, U, \mathbf y, V, t-\tau) C_k(\mathbf x, \mathbf y) \frac{\text d r_\alpha (U)}{\text d U} \text d \mathbf y \text d V \text d \tau  \\
    & \mathcal D_{ij} \approx \delta_{i,d+1} \delta_{j,d+1} r_\alpha(U) \int_0^t \int_{\tilde \Omega} G(\mathbf x, U, \mathbf y, V, t-\tau) C_k(\mathbf x, \mathbf y)   r_\alpha (V) \text d \mathbf y \text d V \text d \tau, \quad i,j = 1, \dots, d + 1. \notag
\end{align}
Here $\delta_{i,d+1}$ is the Kronecker delta, $C_k(\mathbf x, \mathbf y) = \langle k'(\mathbf x') k'(\mathbf x) \rangle$ is the covariance function of $k'(\mathbf x)$, and $G(\mathbf x,U, \mathbf y, V,t-\tau)$ is the the mean-field Green's function that is defined as a solution of
\begin{align}\label{eq:green}
    \frac{\partial G}{\partial \tau} + \mathbf v \cdot \nabla' G + \langle k \rangle \frac{\text d r_\alpha G}{\text d U} = - \delta ( \mathbf x - \mathbf y) \delta ( U - V) \delta (t-\tau), \qquad \tau < t
\end{align} 
with homogeneous initial (at $\tau = 0$) and boundary conditions on $\partial \tilde \Omega$.  
The closure approximations are thus expressed in terms of the mean and two-point covariance of the random input $k(\mathbf x)$.

The derivation of~\eqref{eq:CDF} and~\eqref{eq:closuresG} is based on the following assumptions: $\nabla F$ varies slowly in space and time to justify the use of a local model, the random inputs are mutually uncorrelated, and the variance $\sigma_k^2$ is sufficiently small to warrant its use as a perturbation parameter. 

Our numerical experiments consider one-dimensional ($d=1$) advection in a deterministic velocity field with $v=1$ and linear reaction ($\alpha=1$) with second-order stationary reaction rate $k(\mathbf x)$ with constant mean $\langle k \rangle$ and variance $\sigma_k^2$ and covariance function $C_k(x-y)$. The flow takes place in the semi-infinite domain $\Omega$, so that $\tilde \Omega = [0,\infty) \times [U_{\text{min}}, U_{\text{max}}]$. The deterministic equilibrium state is set to $u_\text{eq}=0$. Under these conditions,~\eqref{eq:closuresG} reduces to
\begin{align}\label{eq:closureterms}
&  \mathcal D_{11} = 0, \quad \mathcal D_{12} = \mathcal D_{21} = 0, \quad \mathcal D_{22} = U^2 \int_0^{t^*} \text{e}^{\langle k \rangle \tau} C_k(v \tau) \text d \tau \notag \\
& \mathcal Q_1 = v, \quad \mathcal Q_2(\mathbf x,U,t) = -U \langle k \rangle + U \int_0^{t^*} \text{e}^{\langle k \rangle \tau} C_k(v\tau)\text d \tau.
\end{align}
where $t^* = \min \{ t, \langle k \rangle^{-1} \log (U_{\text{max}}/U) , x/v \}$. We consider three models of spatial correlation of $k(\mathbf x)$. The first takes $k(\mathbf x)$ to be perfectly correlated, so that $C_k(x-y) = \sigma_k^2$;
then~\eqref{eq:closureterms} simplify to~\eqref{eq:closurek}. The second considers the opposite case, i.e. uncorrelated random field with $C_k(x-y) = \sigma_k^2 \delta (x-y)$, which yields  ~\eqref{eq:closureterms_uncorr}.
Finally, the third one deals with the exponential covariance function $C_k(x-y) = \sigma_k^2 \exp(-|x-y|/\lambda_k)$, where $\lambda_k$ is the correlation length of $k(\mathbf x)$, with closure parameters \eqref{eq:closureterms_finite}.

\begin{figure}[htbp]
\centering\includegraphics[width=\textwidth]{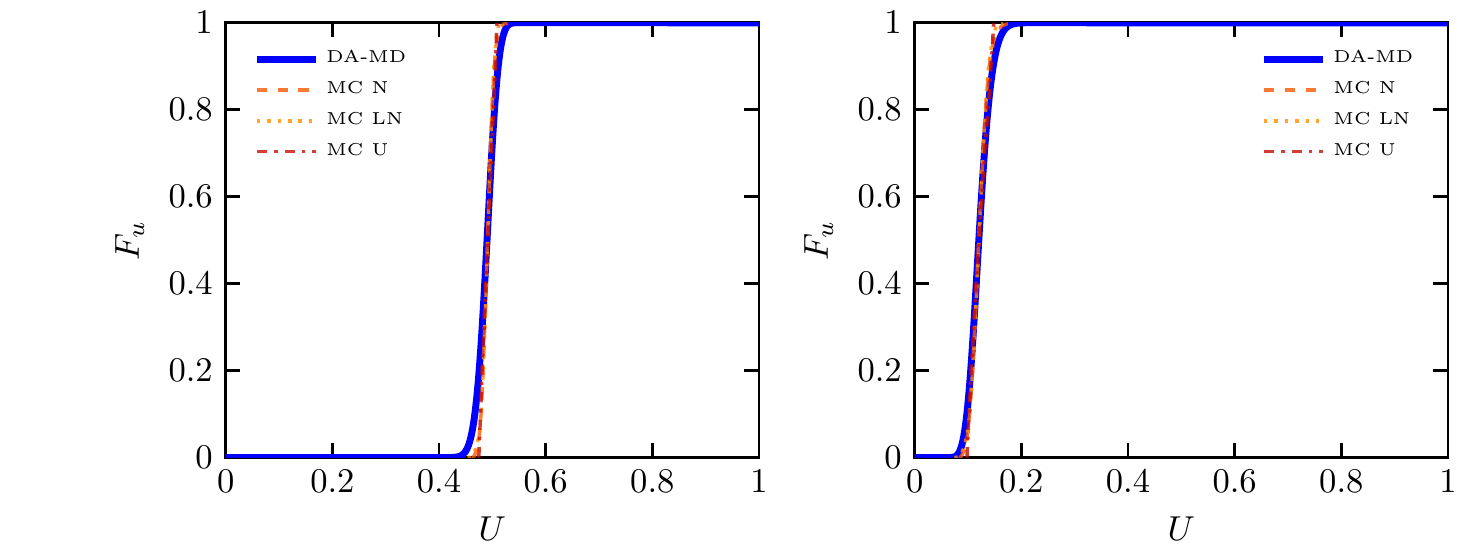}
\caption{Comparison between the FV approximation of the prior CDF and its MC counterpart for the random $k$ scenario. Both techniques use the same mean and variance for $k$, $\langle k \rangle = 2, \sigma_k=0.2$. MC simulations are repeated for different distributions of $k$ sharing the same mean and variance: Normal, Lognormal and Uniform distributions, respectively. Parameters are set to: $N_{\text{MC}} = 1000$, $\Delta t = 0.001$, $\Delta x = 1.6 \cdot 10^{-4}$, $\Delta U = 8.3 \cdot 10^{-4}$. }
\label{fig:CDFprofileMC}
\end{figure}

The CDF equation~\eqref{eq:CDF}, whose coefficients are defined by~\eqref{eq:closureterms}, depends only on the low moments of $k(\mathbf x)$, i.e., on $\langle k \rangle$, $\sigma_k^2$ and $C_k$, rather than on its full PDF.  We study the sensitivity of our closure to a choice of the functional form of the single-point PDF $f_k(K;\mathbf x)$ of $k(\mathbf x)$ for the perfectly correlated case. This is done by comparing a numerical (finite-volume) solution of~\eqref{eq:CDF} with the results of Monte Carlo simulations. The latter consist of post-processing of $N_\text{MC} = 1000$ analytical solutions of the physical model~\eqref{eq:case0}, whose parameters are drawn, alternatively, from the Gaussian, log-normal and uniform PDFs $f_k(K;\mathbf x)$), with negligible discrepancy in CDF terms (see \cref{fig:CDFprofileMC}). As uncertainty is reduced via data assimilation, the discrepancy between posteriors obtained with different assumed PDF forms of $k$ reduces, and the impact of closure approximations on the CDF equation decreases.

\bibliographystyle{unsrt}  


\end{document}